\documentclass{elsart}
\usepackage{amssymb}
\usepackage{amsmath}
\usepackage[english]{babel}

\newtheorem{theorem}{Theorem}

\newtheorem{corollary}[theorem]{Corollary}

\newtheorem{lemma}[theorem]{Lemma}

\newtheorem{remark}[theorem]{Remark}

\begin{document}
\begin{frontmatter}
\title{Semi-classical spectral estimates for Schr\"odinger operators at a
critical level. Case of a degenerate maximum of the potential.}
\author{Brice Camus}
\thanks[label2]{This work was supported by the
\textit{SFB/TR12} project, \textit{Symmetries and Universality in
Mesoscopic Systems} and by the \textit{IHP-Network,
Analysis\&Quantum} reference HPRN-CT-2002-00277. The author thanks
Daniel Barlet for helpful comments on degenerate oscillatory
integrals.}

\address{Ruhr-Universit\"at Bochum, Fakult\"at f\"ur Mathematik,\newline%
Universit\"atsstr. 150, D-44780 Bochum, Germany.\newline%
Email : brice.camus@univ-reims.fr}

\begin{abstract}
We study the semi-classical trace formula at a critical energy
level for a Schr\"odinger operator on $\mathbb{R}^{n}$. We assume
here that the potential has a totally degenerate critical point
associated to a local maximum. The main result, which establishes
the contribution of the associated equilibrium in the trace
formula, is valid for all time in a compact subset of $\mathbb{R}$
and includes the singularity in $t=0$. For these new contributions
the asymptotic expansion involves the logarithm of the parameter
$h$. Depending on an explicit arithmetic condition on the
dimension and the order of the critical point, this logarithmic
contribution can appear in the leading term.
\end{abstract}

\end{frontmatter}

\section{Introduction.}
Let us consider $P_{h}$ a self-adjoint $h$-pseudodifferential
operator, or more generally $h$-admissible (see \cite{[Rob]}),
acting on a dense subset of $L^2(\mathbb{R}^n)$. A classical and
accessible problem is to study the asymptotic behavior, as $h$
tends to 0, of the spectral distributions :
\begin{equation}
\gamma (E,h,\varphi )=\sum\limits_{|\lambda _{j}(h)-E|\leq
\varepsilon }\varphi (\frac{\lambda _{j}(h)-E}{h}),  \label{Def
trace}
\end{equation}
where the $\lambda _{j}(h)$ are the eigenvalues of $P_{h}$, $E$ is
an energy level of the principal symbol of $P_h$ and $\varphi$ a
function. Here we suppose that the spectrum is discrete in
$[E-\varepsilon ,E+\varepsilon ]$, a sufficient condition for this
is given below. If $p_0$ is the principal symbol of $P_{h}$ we
recall that an energy $E$ is regular when $\nabla p_0(x,\xi )\neq
0$ on the energy surface :
\begin{equation}
\Sigma _{E}=\{(x,\xi )\in T^{\star}\mathbb{R}^{n}\text{ }/\text{
}p_0(x,\xi )=E\},
\end{equation}
and critical when it is not regular.

A classical result is the existence of a link between the
asymptotics of (\ref{Def trace}), as $h$ tends to 0, and the
closed trajectories of the Hamiltonian flow of $p_0$ on the energy
surface $\Sigma_{E}$, i.e. :
\begin{equation*}
\lim_{h\rightarrow 0}\gamma (E,h,\varphi )\rightleftharpoons
\{(t,x,\xi)\in\mathbb{R}\times \Sigma _{E} \text{ / }
\Phi_{t}(x,\xi)=(x,\xi)\},
\end{equation*}
where $\Phi_{t}=\mathrm{exp}(tH_{p_0})$ and $H_{p_0}=\partial
_{\xi }p_0.\partial_{x} -\partial _{x}p_0.\partial_{\xi}$. This
duality between spectrum and periodic orbits exists in a lot of
various settings such as in the Selberg trace formula or for the
trace of the wave operator on compact manifolds \cite{D-G}. In the
semi-classical setting this correspondence was initially pointed
out in the physic literature : Gutzwiller \cite{GUT},
Balian\&Bloch \cite{BB}.

For a rigorous mathematical approach, and when $E$ is a regular
energy, a non-exhaustive list of references is Helffer\&Robert
\cite{HR}, Brummelhuis\&Uribe \cite{BU}, Paul\&Uribe \cite{PU},
and more recently Combescure et al. \cite{CRR}, Petkov\&Popov
\cite{P-P}.

Equilibriums are suspected to give special contributions in both
sides of the trace formula. When $E$ is no more a regular value,
the asymptotic behavior of Eq. (\ref{Def trace}) depends on the
nature of the singularities of $p$ on $\Sigma_{E}$ which is too
complicated to be treated in general position. The case of a
non-degenerate critical energy for $p_0$, that is such that the
critical-set $\mathbb{\frak{C}}(p_0) =\{(x,\xi )\in T^{\ast
}\mathbb{R}^{n}\text{ }/\text{ }dp_0(x,\xi )=0\}$ is a compact
$C^{\infty }$ manifold with a Hessian $d^{2}p_0$ transversely
non-degenerate along this manifold, has been investigated first by
Brummelhuis et al. in \cite{BPU}. They treated this question for
quite general operators but for some ''small times'', i.e. it was
assumed that 0 was the only period of the linearized flow in
$\rm{supp}(\hat{\varphi})$. Later, Khuat-Duy in \cite{KhD1} has
obtained the contributions of the non-zero periods of the
linearized flow for $\rm{supp}(\hat{\varphi})$ compact, but for
Schr\"{o}dinger operators with symbol $\xi ^{2}+V(x)$ and a
non-degenerate potential $V$. Our contribution to this subject was
to generalize his result for some more general operators, always
with $\hat{\varphi}$ of compact support and under some geometrical
assumptions on the flow (see \cite{Cam}). In \cite{Cam3} we have
studied the case of a Schr\"odinger operator near a degenerate
minimum of the potential and the objective of the present work is
to investigate the situation near a degenerate maximum which leads
to a totally different asymptotic problem.

After a reformulation, via the theory of Fourier integral
operators, the spectral distribution of Eq. (\ref{Def trace}) can
be expressed in terms of oscillatory integrals whose phases are
related to the classical flow of $p_0$. Moreover, the asymptotic
behavior as $h$ tends to 0 of these oscillatory integrals is
related to closed orbits of this flow. When $(x_0,\xi_0)$ is a
critical point of $p_0$, and hence an equilibrium of the flow, it
is well known that the relation :
\begin{equation}\label{cond.IO.singular}
\mathfrak{F}_t=\mathrm{Ker}(d_{x,\xi}\Phi_{t}(x_0,\xi_0)-\mathrm{Id})\neq
\{ 0\},
\end{equation}
leads to the study of degenerate oscillatory integrals. In the
present work we consider the case of a Schr\"odinger operator :
\begin{equation}
P_h=-h^2\Delta+V(x),
\end{equation}
but a generalization to an $h$-admissible operator (in the sense
of \cite{[Rob]}) of principal symbol $\xi^2+V(x)$ is outlined in
the last section. In particular, we will consider the case of a
potential $V$ with a single and degenerate critical point $x_0$
attached to a local maximum. A typical example is the top of a
polynomial double well. An immediate consequence is that the
symbol admits a unique critical point $(x_0,0)$ on the energy
surface $\{\xi^2+V(x)=V(x_0)\}$ and that the linearized flow at
this point is given by the flow of the free Laplacian. A fortiori,
Eq. (\ref{cond.IO.singular}) is automatically satisfied with :
\begin{gather*}
\mathfrak{F}_t=\{ (\delta u,\delta v)\in
T_{x_0,\xi_0}T^*\mathbb{R}^n \text{ / } \delta v =0 \} \simeq
\mathbb{R}^n, \text{ } t\neq 0,\\
\mathfrak{F}_0=T_{x_0,\xi_0} (T^*\mathbb{R}^n) \simeq
\mathbb{R}^{2n}.
\end{gather*}
In particular, the stationary phase method cannot be applied at
all in a microlocal neighborhood of $t=0$.

The core of the proof lies in establishing suitable local normal
forms for the local phase functions of a Fourier integral operator
approximating the propagator in the semi-classical regime and in a
generalization of the stationary phase formula for these normal
forms. This generalization, based on an analytic representation of
the associated class of oscillatory integrals, is more complicated
than in the case of a local minimum but however allows to compute,
invariantly, the leading term of the related local trace formula.
\section{Hypotheses and main result.}
Let $p(x,\xi)=\xi^2 +V(x)$ where the potential $V$ is smooth on
$\mathbb{R}^n$ and real valued. To this symbol is attached the
$h$-differential operator $P_h=-h^2\Delta+V(x)$ and by a classical
result $P_h$ is essentially autoadjoint, for $h$ small enough, if
$V$ is bounded from below. Moreover, if $E$ is an energy level of
$p$ satisfying :

$(H_{1})$ \textit{There exists } $\varepsilon_{0}>0$ \textit{ such that }%
$p^{-1}([E-\varepsilon_{0},E+\varepsilon_{0}])$\textit{\ is
compact,}

then, by Theorem 3.13 of \cite{[Rob]} the spectrum $\sigma
(P_{h})\cap [E-\varepsilon ,E+\varepsilon ]$ is discrete and
consists in a sequence $\lambda _{1}(h)\leq \lambda _{2}(h)\leq
...\leq \lambda _{j}(h)$ of eigenvalues of finite multiplicities,
if $\varepsilon$ and $h$ are small enough. For example, $(H_1)$ is
certainly satisfied if $V$ goes to infinity at infinity and, more
generally, this is true when $E<\liminf\limits_{\infty} V<\infty$.

We want to study the asymptotic behavior of the spectral
distribution :
\begin{equation}
\gamma (E_{c},h,\varphi)=\sum\limits_{\lambda _{j}(h)\in
[E_{c}-\varepsilon ,E_{c}+\varepsilon ]}\varphi (\frac{\lambda
_{j}(h)-E_{c}}{h}). \label{Objet trace}
\end{equation}
We use the subscript $E_c$ to recall that this energy level is
critical. To avoid any problem of convergence we impose the
condition : \medskip\\
$(H_{2})$\textit{ We have }$\hat{\varphi}\in C_{0}^{\infty
}(\mathbb{R})$ \textit{ with a sufficiently small support near the
origin.}
\begin{remark}\rm{This extra condition on the size of the support is simply here
to avoid contributions of non-trivial closed orbits and can be
easily relaxed. An explicit characterization of
$\mathrm{supp}(\hat{\varphi})$ is given in Lemma \ref{periods}.}
\end{remark}
To simplify notations we write $z=(x,\xi)$ for any point of the
phase space and let be $\Sigma _{E_{C}}=p^{-1}(\{E_{c}\})$. In the
next condition degenerate means that the second derivative at the
critical point $x_0$ is zero. We impose now the type of singularity of the potential :\medskip\\
$(H_{3})$\textit{ On }$\Sigma _{E_{C}}$ \textit{the symbol
}$p$\textit{ has a unique critical point }$z_{0}=(x_{0},0).$
\textit{This critical point is degenerate and associated to a
local maximum of the potential $V$ of the form :
\begin{equation}\label{form pot}
V(x)=E_c+ V_{2k}(x)+\mathcal{O}(||x-x_0||^{2k+1}),
\end{equation}
where $V_{2k}$ is homogeneous of degree $2k$ w.r.t. $(x-x_0)$.
Also, $k\geq 2$ and $V_{2k}$ is definite negative.}
\begin{remark} \label{critical point}
\rm{Since all previous derivatives are 0 in $x_0$, the function
$V_{2k}$ does not depend on the choice of local coordinates near
$x_0$. Contrary to the case exposed in \cite{Cam3}, $z_0$ is not
isolated on $\Sigma_{E_c}$ and might eventually interfere with
non-trivial closed trajectories. A generalization to more general
maximums, e.g. to a sum of such homogeneous terms with different
degrees, is possible but to simplify we only consider the
homogeneous case $(H_3)$. $\hfill{\square}$}
\end{remark}
In this work, we are particulary interested in the contribution of
the fixed point $z_{0}$. To understand the new phenomenon it
suffices to study the localized problem :
\begin{equation} \label{trace-local}
\gamma _{z_{0}}(E_{c},h,\varphi)=\frac{1}{2\pi }\mathrm{Tr}\int\limits_{\mathbb{R}}e^{i%
\frac{tE_{c}}{h}}\hat{\varphi}(t)\psi ^{w}(x,hD_{x})\mathrm{exp}(-\frac{it}{h}%
P_{h})\Theta (P_{h})dt.
\end{equation}
Here $\Theta $ is a function of localization near the critical
energy surface $\Sigma_{E_c}$, $\psi \in C_{0}^{\infty }(T^{\ast
}\mathbb{R}^{n})$ is micro-locally supported near $z_{0}$ and
$\psi ^{w}(x,hD_{x})$ the associated operator obtained by $h$-Weyl
quantization. Rigorous justifications are given in section 3 for
the introduction of $\Theta (P_{h})$ and in section 4 for $\psi
^{w}(x,hD_{x})$. In \cite{Cam3} it was proven that :
\begin{theorem}
Under $(H_{1})$, if $x_0$ is a local minimum of the potential ,
homogeneous as in $(H_3)$, then for all $\varphi$ with
$\hat{\varphi}\in C_0^{\infty}(\mathbb{R})$ we have :
\begin{equation*}
\gamma _{z_{0}}(E_{c},h,\varphi)\sim
h^{-n+\frac{n}{2}+\frac{n}{2k}} \sum\limits_{j,l\in\mathbb{N}^2}
h^{\frac{j}{2}+\frac{l}{2k}}\Lambda _{j,l}(\varphi ),
\end{equation*}
where the $\Lambda _{j,l}$ are some computable distributions. The
leading coefficient is :
\begin{equation}
\Lambda _{0,0}(\varphi )
=\frac{\mathrm{S}(\mathbb{S}^{n-1})}{(2\pi)^n}
\int\limits_{\mathbb{S}^{n-1}} |V_{2k}(\eta)|^{-\frac{n}{2k}}
d\eta \int\limits_{\mathbb{R}_{+} \times \mathbb{R}_{+}}
\varphi(u^2 +v^{2k}) u^{n-1} v^{n-1} dudv,
\end{equation}
where $\mathrm{S}(\mathbb{S}^{n-1})$ is the surface of the
unit-sphere of $\mathbb{R}^n$.
\end{theorem}
Clearly, this result shows that the spectral estimates are related
to a purely local problem, namely to estimate the asymptotic
behavior of fiber integrals :
\begin{equation}
\int\limits_{T^*\mathbb{R}^n} \varphi (\frac{\xi^2
+V_{2k}(x)}{h})dxd\xi,\text{ }h\rightarrow 0^{+} .
\end{equation}
This result can be interpreted as a scaling of the trace of
$-\Delta+V_{2k}$, cf. the first term of the trace formula.
Evidently, such an interpretation cannot hold in our setting since
the trace does not exists if $V_{2k}$ is a negative definite
function.

The main result of the present work is :
\begin{theorem}\label{Main}
Under hypotheses $(H_{1})$ to $(H_{3})$ we have :
\begin{equation*}
\gamma _{z_{0}}(E_{c},h,\varphi)\sim \sum\limits_{m=0,1}
h^{-n+\frac{n}{2}+\frac{n}{2k}} \sum\limits_{j,l\in\mathbb{N}^2}
h^{\frac{j}{2}+\frac{l}{k}}\mathrm{log}(h)^m \Lambda
_{j,l,m}(\varphi ),
\end{equation*}
where the $\Lambda _{j,l,m}$ are in $\mathcal{S}'(\mathbb{R})$.\\
As concerns the leading term of the expansion, when
$\frac{n(k+1)}{2k}\notin \mathbb{N}$, the first non-zero
coefficient of this local trace formula is given by :
\begin{equation}
h^{-n+\frac{n}{2} +\frac{n}{2k}}%
\left\langle T_{n,k},\varphi\right\rangle%
\frac{\mathrm{S}(\mathbb{S}^{n-1})}{(2\pi)^n}
\int\limits_{\mathbb{S}^{n-1}} |V_{2k}(\eta)|^{-\frac{n}{2k}}d\eta
.
\end{equation}
The distributions $T_{n,k}$ are given by :
\begin{gather}
\left\langle T_{n,k},\varphi\right\rangle
=\int\limits_{\mathbb{R}}
(C_{n,k}^{+}|t|_{+}^{n\frac{k+1}{2k}-1}+C_{n,k}^{-}|t|_{-}^{n\frac{k+1}{2k}-1}) \varphi (t) dt,
\text{ if }n \text{ is odd},\\
\left\langle T_{n,k},\varphi\right\rangle
=C_{n,k}^{-}\int\limits_{\mathbb{R}} |t|_{-}^{n\frac{k+1}{2k}-1}
\varphi (t) dt,\text{ if }n \text{ is even}.
\end{gather}
But if $\frac{n(k+1)}{2k}\in \mathbb{N}$ and $n$ is odd then the
top-order term is :
\begin{equation}
C_{n,k} \log (h)h^{-n+\frac{n}{2}
+\frac{n}{2k}}\frac{\mathrm{S}(\mathbb{S}^{n-1})}{(2\pi)^n}\int\limits_{\mathbb{S}^{n-1}}
|V_{2k}(\eta)|^{-\frac{n}{2k}} d\eta  \int\limits_{\mathbb{R}}
|t|^{n\frac{k+1}{2k}-1} \varphi (t) dt.
\end{equation}
Finally, if $\frac{n(k+1)}{2k}\in \mathbb{N}$ and $n$ is even,
$C_{n,k}^{+}=C_{n,k}^{-}$ and we have :
\begin{equation}
C_{n,k}^{\pm} h^{-n+\frac{n}{2} +\frac{n}{2k}}\frac{1}{(2\pi)^n}
\int\limits_{\mathbb{S}^{n-1}} |V_{2k}(\eta)|^{-\frac{n}{2k}}d\eta
\int\limits_{\mathbb{R}} |t|^{n\frac{k+1}{2k}-1} \varphi (t) dt.
\end{equation}
In all expressions above $C_{n,k}$, $C_{n,k}^{\pm}$ are non-zero
universal constants depending only on $n$ and $k$.
\end{theorem}
Explicit formulas for the numbers $C^{\pm}_{n,k}$ are given in
section 6. The arithmetical condition on $k$ and $n$ might be
surprising at the first look. But this condition becomes clear
when the oscillatory integrals of our spectral problem are
analytically reformulated in section 6. Also, viewing the top
order coefficient of the trace as a tempered distribution on the
Schwartz function $\varphi$, we obtain that the singular support
is only located at the origin.
\section{Oscillatory representation of the spectral functions.}
The construction below is more or less classical and will be
sketchy. For a more detailed exposition the reader can consult
\cite{BPU}, \cite{Cam} or \cite{KhD1}.

Let be $\varphi \in \mathcal{S}(\mathbb{R})$ with
$\hat{\varphi}\in C_{0}^{\infty }(\mathbb{R})$, we recall that :
\begin{equation*} \gamma (E_{c},h,\varphi)
=\sum\limits_{\lambda _{j}(h)\in I_{\varepsilon }}\varphi
(\frac{\lambda _{j}(h)-E_{c}}{h}),
\end{equation*}
where $I_{\varepsilon } =[E_{c}-\varepsilon ,E_{c}+\varepsilon ]$,
with $0<\varepsilon<\varepsilon_0$, and
$p^{-1}(I_{\varepsilon_{0}})$ compact in $T^{\ast
}\mathbb{R}^{n}$. We localize near the critical energy level
$E_{c}$ by inserting a cut-off function $\Theta \in C_{0}^{\infty
}(]E_{c}-\varepsilon ,E_{c}+\varepsilon \lbrack )$, such that
$\Theta =1$ in a neighborhood of $E_{c}$ and $0\leq \Theta \leq 1$
on $\mathbb{R}$. The corresponding decomposition is :
\begin{equation*}
\gamma (E_{c},h,\varphi)=\gamma _{1}(E_{c},h,\varphi)+\gamma
_{2}(E_{c},h,\varphi),
\end{equation*}
with :
\begin{equation}
\gamma _{1}(E_{c},h,\varphi)=\sum\limits_{\lambda _{j}(h)\in
I_{\varepsilon }}(1-\Theta )(\lambda _{j}(h))\varphi
(\frac{\lambda _{j}(h)-E_{c}}{h}),
\end{equation}
\begin{equation}
\gamma _{2}(E_{c},h,\varphi)=\sum\limits_{\lambda _{j}(h)\in
I_{\varepsilon }}\Theta (\lambda _{j}(h))\varphi (\frac{\lambda
_{j}(h)-E_{c}}{h}).
\end{equation}
Since $\varphi\in \mathcal{S}(\mathbb{R})$ a classical estimate,
see e.g. Lemma 1 of \cite{Cam1}, is :
\begin{equation}
 \gamma _{1}(E_{c},h,\varphi)=\mathcal{O}(h^{\infty }),
\text{ as } h\rightarrow 0.\label{S1(h)=Tr}
\end{equation}
By inversion of the Fourier transform we have :
\begin{equation*}
\Theta (P_{h})\varphi (\frac{P_{h}-E_{c}}{h})=\frac{1}{2\pi}\int\limits_{%
\mathbb{R}}e^{i\frac{tE_{c}}{h}}\hat{\varphi}(t)\mathrm{exp}(-\frac{it}{h}%
P_{h})\Theta (P_{h})dt.
\end{equation*}
The trace of the left hand-side is precisely $\gamma
_{2}(E_{c},h,\varphi)$ and Eq. (\ref{S1(h)=Tr}) gives :
\begin{equation}\label{Trace S2(h)}
\gamma (E_{c},h,\varphi)=\frac{1}{2\pi }\mathrm{Tr}\int\limits_{\mathbb{R}}e^{i%
\frac{tE_{c}}{h}}\hat{\varphi}(t)\mathrm{exp}(-\frac{it}{h}P_{h})\Theta
(P_{h})dt+\mathcal{O}(h^{\infty }).
\end{equation}
If $U_{h}(t)=\mathrm{exp}(-\frac{it}{h}P_{h})$ is the evolution
operator, we can approximate $U_{h}(t)\Theta (P_{h})$ by a Fourier
integral-operator depending on a parameter $h$. If $\Lambda$ is
the Lagrangian manifold associated to the flow of $p$ :
\begin{equation*}
\Lambda =\{(t,\tau ,x,\xi ,y,\eta )\in T^{\ast }\mathbb{R}\times T^{\ast }%
\mathbb{R}^{n}\times T^{\ast }\mathbb{R}^{n}:\tau =p(x,\xi
),\text{ }(x,\xi )=\Phi _{t}(y,\eta )\},
\end{equation*}
a classical and general result, see e.g. Duistermaat \cite{DUI1},
is :
\begin{theorem}
The operator $U_{h}(t)\Theta (P_{h})$ is an $h$-FIO associated to
$\Lambda$. For each $N\in\mathbb{N}$ there exists $U_{\Theta
,h}^{(N)}(t)$ with integral kernel in H\"ormander's class
$I(\mathbb{R}^{2n+1},\Lambda )$ and $R_{h}^{(N)}(t)$ bounded, with
a $L^{2}$-norm uniformly bounded for $0<h\leq 1$ and $t$ in a
compact subset of $\mathbb{R}$, such that :
\begin{equation*}
U_{h}(t)\Theta(P_{h})=U_{\Theta,h}^{(N)}(t)+h^{N}R_{h}^{(N)}(t).
\end{equation*}
\end{theorem}
The remainder associated to $R_{h}^{(N)}(t)$ is controlled by the
classical trick :
\begin{corollary}
Let $\Theta _{1}\in C_{0}^{\infty }(\mathbb{R})$, with $\Theta
_{1}=1$ on $\rm{supp}(\Theta )$ and $\rm{supp}(\Theta _{1})\subset
I_{\varepsilon }$, then $\forall N\in \mathbb{N}$ :
\begin{equation*}
\mathrm{Tr}(\Theta (P_{h})\varphi (\frac{P_{h}-E_{c}}{h}))=\frac{1}{2\pi }%
\mathrm{Tr}\int\limits_{\mathbb{R}}\hat{\varphi}(t)e^{\frac{i}{h}%
tE_{c}}U_{\Theta ,h}^{(N)}(t)\Theta
_{1}(P_{h})dt+\mathcal{O}(h^{N-n}).
\end{equation*}
\end{corollary}
The proof is easy by cyclicity of the trace (see \cite{Cam1} or \cite{[Rob]}).\medskip\\
For the particular case of a Schr\"odinger operator the BKW ansatz
shows that the integral kernel of $U_{\Theta,h}^{(N)}(t)$ can be
recursively constructed as :
\begin{gather*}
K_h^{(N)}(t,x,y)=\frac{1}{(2\pi h)^n} \int\limits_{\mathbb{R}^n}
b_h^{(N)}(t,x,y,\xi) e^{\frac{i}{h} (S(t,x,\xi)-\left\langle y,\xi
\right\rangle)} d\xi,\\
b_h^{(N)}=b_0+h b_1+...+h^N b_N,
\end{gather*}
where $S$ satisfies the Hamilton-Jacobi equation :
\begin{equation*}
\left\{
\begin{array}{c}
\partial _{t}S(t,x,\xi )+ p(x,\partial
_{x}S(t,x,\xi))=0, \\
S(0,x,\xi)=\left\langle x,\xi \right\rangle .
\end{array}
\right.
\end{equation*}
In particular we obtain that :
\begin{equation*}
\{(t,\partial _{t}S(t,x,\eta ),x,\partial _{x}S(t,x,\eta
),\partial _{\eta }S(t,x,\eta ),-\eta )\}\subset
\Lambda_{\mathrm{flow}} ,
\end{equation*}
and that the function $S$ is a generating function of the flow,
i.e. :
\begin{equation}
\Phi _{t}(\partial _{\eta }S(t,x,\eta ),\eta ) =(x,\partial
_{x}S(t,x,\eta )). \label{Gene}
\end{equation}
Inserting this approximation in Eq. (\ref{Trace S2(h)}) we find
that, modulo an error $\mathcal{O}(h^{N-n})$, the trace
$\gamma(E_{c},h,\varphi)$ can be written for all $N\in \mathbb{N}$
as :
\begin{equation}
\gamma(E_{c},h,\varphi)=\sum\limits_{j<N}\frac{h^j}{(2\pi h)^n}
\int\limits_{\mathbb{R}\times
T^*\mathbb{R}^n}e^{\frac{i}{h}(S(t,x,\xi )-\left\langle x,\xi
\right\rangle +tE_{c})}a_{j}(t,x,\xi )\hat{\varphi}(t)dtdxd\xi,
\label{gamma1 OIF}
\end{equation}
where $a_{j}(t,x,\eta )=b_{j}(t,x,x,\eta )$ is the evaluation of
$b_j$ on the diagonal $\{x=y\}$.
\begin{remark}
\rm{By a theorem of Helffer\&Robert, see \cite{[Rob]} Theorem 3.11
and Remark 3.14, $\Theta (P_{h})$ is $h$-admissible with a symbol
supported in $p^{-1}(I_{\varepsilon })$. This allows to consider
only oscillatory-integrals with compact support.$\hfill{\square}$}
\end{remark}
\section{Classical dynamic near the equilibrium.}
Critical points of the phase function of (\ref{gamma1 OIF}) are
given by the equations :
\begin{equation*}
\left\{
\begin{array}{c}
E_{c}=-\partial _{t}S(t,x,\xi ), \\
x=\partial _{\xi }S(t,x,\xi ), \\
\xi =\partial _{x}S(t,x,\xi ),
\end{array}
\right. \Leftrightarrow \left\{
\begin{array}{c}
p(x,\xi )=E_{c}, \\
\Phi _{t}(x,\xi )=(x,\xi ),
\end{array}
\right.
\end{equation*}
where the right hand side defines a closed trajectory of the flow
inside $\Sigma _{E_{c}}$. Since we are mainly interested in the
contribution of the critical point, we choose a function $\psi \in
C_{0}^{\infty }(T^{\ast }\mathbb{R}^{n})$, with $\psi =1\text{
near }z_{0}$, hence :
\begin{gather*}
\gamma _{2}(E_{c},h,\varphi) =\frac{1}{2\pi }\mathrm{Tr}\int\limits_{\mathbb{R}}e^{i%
\frac{tE_{c}}{h}}\hat{\varphi}(t)\psi ^{w}(x,hD_{x})\mathrm{exp}(-\frac{i}{h}%
tP_{h})\Theta (P_{h})dt\\
+\frac{1}{2\pi }\mathrm{Tr}\int\limits_{\mathbb{R}}e^{i\frac{tE_{c}}{h}}\hat{%
\varphi}(t)(1-\psi
^{w}(x,hD_{x}))\mathrm{exp}(-\frac{i}{h}tP_{h})\Theta (P_{h})dt.
\end{gather*}
Under the additional hypothesis of having a clean flow, the
asymptotics of the second term is given by the semi-classical
trace formula on a regular level. We also observe that the
contribution of the first term, which is precisely the
distribution $\gamma_{z_0}(E_c,h,\varphi)$ of Theorem \ref{Main},
is micro-local. Hence this allows to introduce local coordinates
near $z_{0}$. As pointed out in Remark \ref{critical point},
$\{z_{0}\}$ is not a connected component of the energy surface
$\Sigma_{E_c}$ and elements :
\begin{equation}
\{(T,z), T\neq 0,\text{ } z\in\Sigma_{E_c}\cap
\mathrm{supp}(\psi),\text{ } \Phi_T(z)=z \},
\end{equation}
could contribute in the asymptotic expansion. In fact, the only
contribution arises from the set $\{(t,z_0),\text{
}t\in\mathrm{supp}(\hat{\varphi})\}$ if the support of
$\hat{\varphi}$ is small enough.
\begin{lemma}\label{periods} There exists a $T>0$, depending only on $V$, such that $\Phi_t(z)\neq z$ for all $z\in
\Sigma_{E_c}\backslash \{z_0\}$ and all $t\in]-T,0[\cup ]0,T[$.
\end{lemma}
\textit{Proof.} If $H_p$ is our hamiltonian vector field and
$z=(x,\xi)$ we have :
\begin{equation*}
|| H_p(z_1)-H_p(z_2)||^2= 4||\xi_1-\xi_2||^2
+||\nabla_xV(x_1)-\nabla_xV(x_2)||^2.
\end{equation*}
Since our potential is smooth, when $z_1$ and $z_2$ are in in the
energy surface $\Sigma_{E_c}$, which is compact by assumption,
there exists $M>0$ such that :
\begin{equation*}
||H_p(z_1)-H_p(z_2)||\leq  M ||z_1-z_2||.
\end{equation*}
The main result of \cite{Yor} shows that any periodic trajectory
inside $\Sigma_{E_c}$ has a period $p$ such that $p\geq
\frac{2\pi}{M}$. $\hfill{\blacksquare}$
\begin{remark} \rm{The result of Lemma \ref{periods} can perhaps be improved for any $T>0$ by replacing
$\Sigma_{E_c}$ by $U_T\cap \Sigma_{E_c}$, where $U_T$ is a
restricted neighborhood of $z_0$. This is suggested by the fact
that such a result is proven in \cite{KhD1} for non-degenerate
potentials. Anyhow, via Lemma \ref{periods} we can understand the
contribution in $t=0$ of elements $\{t,z_0\}$ in the trace formula
if $\hat{\varphi}\in C_0^{\infty}(]-p,p[)$.$\hfill{\square}$}
\end{remark}
Now, we restrict our attention to the singular contribution
generated by the critical point. Since $z_0$ is an equilibrium of
the flow we obtain that :
\begin{equation}
d_{x,\xi}\Phi _{t}(z_{0})=\mathrm{exp}(tH_{-\Delta}),\text{
}\forall t.
\end{equation}
The computation of this linear operator is easy and gives :
\begin{equation}
d_{x,\xi}\Phi_{t}(z_{0})(u,v)=(u+2tv,v),\text { } \forall (u,v)\in
T_{z_0} T^*\mathbb{R}^n.
\end{equation}
From classical mechanics we know that the singularities, in the
sense of the Morse theory, of the function
$S(t,x,\xi)-\left\langle x,\xi\right\rangle$ are supported in the
set :
\begin{equation*}
\mathfrak{F}_t=\mathrm{Ker}(d_{x,\xi}\Phi_t(z_0)-\mathrm{Id}),
\end{equation*}
see e.g. Lemma 9 of \cite{Cam}. As mentioned in the introduction
we obtain :
\begin{gather*}
\mathfrak{F}_0 = T_{z_0} (T^* \mathbb{R}^{n})\simeq \mathbb{R}^{2n},\\
\mathfrak{F}_t = \{ (u,v) \in T_{z_0} (T^* \mathbb{R}^{n}) \text{
/ } v=0 \} \simeq \mathbb{R}^{n}, \text{ }t\neq 0.
\end{gather*}
To simplify notations, and until further notice, all derivatives
will be taken with respect to the initial conditions $(x,\xi)$.
The next non-zero terms of the Taylor expansion of the flow are
computed via the technical result :
\begin{lemma}
\label{TheoFormule de récurence du flot}Let be $z_{0}$ an
equilibrium of the $C^{\infty}$ vector field $X$ and $\Phi _{t}$
the flow of $X$. Then for all $m\in \mathbb{N}^{\ast }$, there
exists a polynomial map $P_{m}$, vector valued and of degree at
most $m$, such that :
\begin{equation*}
d^{m}\Phi _{t}(z_{0})(z^{m})=d\Phi
_{t}(z_{0})\int\limits_{0}^{t}d\Phi _{-s}(z_{0})P_{m}(d\Phi
_{s}(z_{0})(z),...,d^{m-1}\Phi _{s}(z_{0})(z^{m-1}))ds.
\end{equation*}
\end{lemma}
For a proof we refer to \cite{Cam} or \cite{Cam1}. Here we shall
use that our vector field is :
\begin{equation*}
H_p= 2\xi \frac{\partial}{\partial x} - \partial_x
V(x)\frac{\partial}{\partial \xi}.
\end{equation*}
We identify the linearized flow in $z_0$ with a matrix
multiplication operator :
\begin{equation*}
d\Phi_t (z_0)=
\begin{pmatrix}
1 & 2t \\
0 & 1
\end{pmatrix}.
\end{equation*}
Clearly, with the hypothesis $(H_3)$ we obtain the polynomials :
\begin{gather*}
P_j=0,\text{ }\forall j\in\{2,..,2k-2\},\\
P_{2k-1}(Y_1,...,Y_{2k-2})=
\begin{pmatrix}
0\\
d^{2k-1} \nabla_x V (x_0)(Y_1 ^{2k-1})
\end{pmatrix}
\neq 0.
\end{gather*}
Where the notation $Y_1^{l}$ stands for $(Y_1,...,Y_1)$ :
$l$-times. Inserting the definition of $d\Phi_s(z_0)$ and
integration from 0 to $t$ yields :
\begin{equation*}
d^{2k-1} \Phi_t(z_0)((x,\xi)^{2k-1})=\begin{pmatrix}
1 & 2t \\
0 & 1
\end{pmatrix}
\int\limits_{0}^{t}
\begin{pmatrix}
2s d^{2k-1} \nabla_x V (x_0)((x+2s\xi)^{2k-1})\\
-d^{2k-1} \nabla_x V (x_0)((x+2s\xi)^{2k-1})
\end{pmatrix}
ds.
\end{equation*}
Terms of higher degree can be obtained similarly by successive
integrations. If we assume that $z_0=0$ the jet of order $2k-1$ of
the flow is :
\begin{equation}\label{germ-flot}
\Phi_t(z)=d\Phi_t(0)(z)+\frac{1}{(2k-1)!}
d^{2k-1}\Phi_t(0)(z^{2k-1})+\mathcal{O}(||z||^{2k}),
\end{equation}
and can be computed explicitly with a given $V_{2k}$.
\section{Normal forms of the phase function.}
Since the contribution we study is local, cf. the introduction of
$\psi$ in Eq. (\ref{trace-local}), we work with some local
coordinates $(x,\xi)$ near the critical point $z_0$. With these
coordinates we identify locally $T^{\ast }\mathbb{R}^{n}\cap
V(z_0)$ with an open of $\mathbb{R}^{2n}$. We define :
\begin{equation}
\Psi(t,z)=\Psi(t,x,\xi)=S(t,x,\xi )-\left\langle x,\xi
\right\rangle+tE_{c},\text{ }z=(x,\xi)\in \mathbb{R}^{2n}.
\label{defphase}
\end{equation}
We start by a more precise description of our phase function.
\begin{lemma}\label{structure phase} Near $z_0$, here supposed to be 0 to simplify, we
have :
\begin{equation}
\Psi(t,x,\xi)=-t||\xi||^2+ S_{2k}(t,x,\xi)+R_{2k+1}(t,x,\xi ),
\label{forme phase}
\end{equation}
where $S_{2k}$ is homogeneous of degree $2k$ w.r.t. $(x,\xi)$ and
is uniquely determined by $V_{2k}$. Moreover,
$R_{2k+1}(t,x,\xi)=\mathcal{O}(||(x,\xi )||^{2k+1})$, uniformly
for $t$ in a compact subset of $\mathbb{R}$.
\end{lemma}
\textit{Proof.} With the particular structure of the flow in
$z_0$, cf. Eq. (\ref{germ-flot}), we search our local generating
function as :
\begin{equation*}
S(t,x,\xi)=-tE_c+S_2(t,x,\xi)+S_{2k}(t,x,\xi)+\mathcal{O}(||(x,\xi)||^{2k+1}),
\end{equation*}
where the $S_j$ are time dependant and homogeneous of degree $j$
w.r.t. $(x,\xi)$. Starting from the implicit relation
$\Phi_t(\partial_{\xi}S(t,x,\xi),\xi)=(x,\partial_{x}S(t,x,\xi))$
and with Eq. (\ref{germ-flot}), we obtain that :
\begin{equation*}
S_2(t,x,\xi)=\left\langle x,\xi \right\rangle-t||\xi ||^2.
\end{equation*}
Now we compute the term $S_{2k}$. To do so, we retain only terms
homogeneous of degree $2k-1$ and we get :
\begin{equation*}
d\Phi_t(0)((\partial_\xi S,0))+\frac{1}{(2k-1)!}d^{2k-1}
\Phi_t(0)((\partial_\xi S_2,\xi)^{2k-1})=(0,\partial_x S_{2k}).
\end{equation*}
If $J$ is the matrix of the usual symplectic form $\sigma$ on
$\mathbb{R}^{2n}$, we have :
\begin{equation*}
J\nabla S_{2k}(t,x,\xi)=\frac{1}{(2k-1)!}d^{2k-1}\Phi
_{t}(0)((x-2t\xi,\xi)^{2k-1}).
\end{equation*}
By homogeneity and with Eq. (\ref{germ-flot}) we obtain :
\begin{equation*}
S_{2k}(t,x,\xi)=\frac{1}{(2k)!}\sigma ((x,\xi),d^{2k-1} \Phi_t(0)
((x-2t\xi,\xi)^{2k-1})).
\end{equation*}
This gives the result since $d^{2k-1} \Phi_t(0)$ is well
determined
by $V_{2k}$. $\hfill{\blacksquare}$\medskip\\
Fortunately, we will not have to compute the remainder explicitly
because of some homogenous considerations. To prepare the
construction of our normal forms, we study carefully $S_{2k}$ .
\begin{corollary} \label{structure fine}
The function $S_{2k}$ satisfies :
\begin{equation} \label{series S_2k}
S_{2k}(t,x,\xi)=-t V_{2k}(x) +t^2 \left\langle \xi, \nabla_x
V_{2k}(x)\right\rangle +\sum\limits_{j,l=1}^{n} \xi_j\xi_l
g_{j,l}(t,x,\xi),
\end{equation}
where the functions $g_{jl}$ are smooth and vanish in $x=0$.
\end{corollary}
\textit{Proof.} If $f$ is homogeneous of degree $2k>2$ w.r.t.
$(x,\xi)$ we can write :
\begin{equation*}
f(t,x,\xi)=f_1(t,x)+\left\langle \xi, f_2(t,x) \right\rangle +
\sum\limits_{j,l=1}^{n} \xi_j\xi_l f_3^{(jl)}(t,x).
\end{equation*}
It remains to compute the function $f_1$ and the vector field
$f_2$. We have :
\begin{equation*}
S_{2k}(t,x,\xi)=\frac{1}{(2k)!}
\sigma((x,\xi),d\Phi_t(0)\int\limits_{0}^{t}
\begin{pmatrix}
2s d^{2k-1} \nabla_x V(0)((x+2(s-t)\xi)^{2k-1})\\
-d^{2k-1} \nabla_x V(0)((x+2(s-t)\xi)^{2k-1})
\end{pmatrix}
ds.
\end{equation*}
Clearly, the term of degree homogeneous of degree $2k$ w.r.t. $x$
is given by :
\begin{equation*}
-\frac{1}{(2k)!}\int\limits_{0}^t \left\langle x, d^{2k-1}
\nabla_x V(0)(x^{2k-1})\right\rangle ds=-\frac{t}{2k} \left\langle
x, \nabla_xV_{2k}(x) \right\rangle =-tV_{2k}(x).
\end{equation*}
Where the last result holds by homogeneity. As concerns the linear
term w.r.t. $\xi$, by combinatoric and linear operations, this one
can be written :
\begin{gather*}
\frac{t^2}{(2k)!}\left( (2k-1) \left\langle x,d^{2k-1} \nabla_x
V(0)((x^{2k-2},\xi))\right\rangle +\left\langle \xi, d^{2k-1}
\nabla_x V(0)(x^{2k-1}) \right\rangle\right )\\
=\frac{t^2}{2k}\left( (2k-1) \left\langle \nabla_x V_{2k}(x),\xi
\right\rangle +\left\langle \xi,\nabla_x V_{2k}(x)
\right\rangle\right )=t^2 \left\langle \xi, \nabla_x V_{2k} (x)
\right\rangle.
\end{gather*}
This complete the proof. $\hfill{\blacksquare}$\medskip\\
These two terms can be also derived by a heuristic method.
Starting from the Hamilton-Jacobi equation we obtain
$S(0,x,\xi)=\left\langle x,\xi\right\rangle$,
$\partial_tS(0,x,\xi)=-p(x,\xi)$ and also :
\begin{equation*}
\partial^2_{t,t} S(t,x,\xi) =-\left\langle \partial_\xi p(x,\partial_x S(t,x,\xi)),
\partial^2_{t,x}S(t,x,\xi)\right\rangle.
\end{equation*}
But, for our flow, in $t=0$ we have simply :
\begin{equation*}\partial^2_{t,t} S(0,x,\xi)= 2\left\langle \xi, \partial_x
V(x)\right\rangle.
\end{equation*}
Hence the Taylor expansion in $t=0$ provides a good result with
few calculations. Unfortunately, this approach gives no
information about the degree, w.r.t. $(x,\xi)$, of the remainders
$\mathcal{O}(t^d)$ for each $d\geq 3$.

We have now enough material to establish the normal of our phase
function.
\begin{lemma}\label{FN1}
In a micro-local neighborhood of $z=z_0$, there exists local
coordinates $\chi$ such that :
\begin{equation}
\Psi(t,z) \simeq -\chi_{0}(\chi_1^2-\chi_{2}^{2k}).
\end{equation}
\end{lemma}
\noindent\textit{Proof.} We can here assume that $z_{0}$ is the
origin. We proceed in two steps. First we want to eliminate terms
of high degree. Starting form Eq. (\ref{series S_2k}) we define :
\begin{equation*}
E(t,x,\xi)=\sum\limits_{j,l}\xi_j \xi_l g_{j,l}(t,x,\xi).
\end{equation*}
Since $S(t,x,\xi)-\left\langle x,\xi\right\rangle=\mathcal{O}(t)$,
we have $\Psi(t,z)=\mathcal{O}(t)$. Hence, all terms of the
expansion are $\mathcal{O}(t)$ and this allows to write
$E=t\tilde{E}$. Similarly, we write $R_{2k+1}=t\tilde{R}_{2k+1}$
in Eq. (\ref{forme phase}). To obtain a blow-up of the
singularity, we use polar coordinates $x=r\theta$, $\xi=q\eta$,
$\theta,\eta \in \mathbb{S}^{n-1}(\mathbb{R})$,
$q,r\in\mathbb{R}_{+}$. This induces naturally a jacobian
$r^{n-1}q^{n-1}$. By construction, there exists a function $F$
vanishing in $(q,r)=(0,0)$ such that :
\begin{equation}
\tilde{E}(t,r\theta,q\eta)=q^2 F(t,r,\theta,q,\eta).
\end{equation}
With Lemma \ref{structure phase} and Corollary \ref{structure
fine}, near $z_0$ the phase $\Psi(t,z)$ can be written :
\begin{equation*}
-t\left(q^{2} +r^{2k} V_{2k}(\theta)-t qr^{2k-1}\left\langle
\eta,\nabla V_{2k}(\theta) \right\rangle +q^2
F(t,r,\theta,q,\eta)+\tilde{R}_{2k+1}(t,r\theta,q\eta)\right ).
\end{equation*}
Thanks to the Taylor formula, the remainder $\tilde{R}_{2k+1}$ can
be written as :
\begin{equation*}
\tilde{R}_{2k+1}(t,r\theta,q\eta)=q^2
R_1(t,r,\theta,q,\eta)+r^{2k}R_2(t,r,\theta,q,\eta),
\end{equation*}
where $R_1$ vanishes in $r=0$ and $R_2$ vanishes in $q=0$. We
obtain :
\begin{equation*}
\Psi(t,z) \simeq -t(q^{2}\alpha_1(t,r,\theta,q,\eta) -r^{2k}
\alpha_2(t,r,\theta,q,\eta))+t^2 qr^{2k-1}\left\langle \eta,\nabla
V_{2k}(\theta) \right\rangle,
\end{equation*}
where we have defined :
\begin{gather*}
\alpha_1(t,r,\theta,q,\eta)=(1+R_1+F)(t,r,\theta,q,\eta),\\
\alpha_2(t,r,\theta,q,\eta)=|V_{2k}(\theta)|+R_2(t,r,\theta,q,\eta).
\end{gather*}
Since $|V_{2k}(\theta)|>0$ on $\mathbb{S}^{n-1}$, we can eliminate
$\alpha_1$ and $\alpha_2$ by a local change of coordinates :
\begin{equation} \label{alpha}
(q \alpha_1^{\frac{1}{2}}, r \alpha_2^{\frac{1}{2k}})\rightarrow
(Q,R),
\end{equation}
near $(q,r)=(0,0)$. This is acceptable since the corresponding
Jacobian is :
\begin{equation}\label{jacobian V_2K}
|J(Q,R)|(0,0)=|V_{2k}(\theta)|^{\frac{1}{2k}}\neq 0.
\end{equation}
In these local coordinates, and still using $(r,q)$ instead of
$(R,Q)$, we obtain :
\begin{equation*}
\Psi(t,z)\simeq -t(q^{2} -r^{2k})+t^2
qr^{2k-1}\varepsilon(t,r,\theta,q,\eta),
\end{equation*}
where :
\begin{equation*}
\varepsilon (t,r,\theta,q,\eta)= \left\langle \eta,\nabla
V_{2k}(\theta)
 \right\rangle(\alpha_1^{-\frac{1}{2}}
\alpha_2^{\frac{1-2k}{2k}})(t,r,\theta,q,\eta).
\end{equation*}
In a second time, we eliminate the nonlinear term in $t$. To do
so, we write :
\begin{equation*}
-t(q^2-tqr^{2k-1}\varepsilon)=-t(q-\frac{t}{2}r^{2k-1}\varepsilon)^2
+\frac{1}{4}t^3 r^{4k-2}\varepsilon^2.
\end{equation*}
Now, we can factor out the last term if we use :
\begin{equation*}
\alpha_3(t,r,\theta,q,\eta)=(1-\frac{t^2}{4}r^{4k-2}\varepsilon^2(t,r,\theta,q,\eta)).
\end{equation*}
Then the change of variables $r\rightarrow
\alpha_3^{-\frac{1}{2k}}r$ gives :
\begin{equation*}
\Psi(t,z)\simeq
-t((q-\frac{t}{2}r^{2k-1}\tilde{\varepsilon}(t,r,\theta,q,\eta))^2
-r^{2k}),
\end{equation*}
where :
\begin{equation*}
\tilde{\varepsilon}(t,r,\theta,q,\eta)=
\alpha_3^{\frac{1-2k}{2k}}\varepsilon(t,\alpha_3^{-\frac{1}{2k}}r,t,\theta,q,\eta).
\end{equation*}
Finally, if we define :
\begin{gather}(\chi_0,\chi_2,\chi_3,...,\chi_{2n})(t,r,\theta,q,\eta) =(t,r,\theta,\eta),\\
\chi_1(t,r,\theta,q,\eta)=
q-\frac{t}{2}r^{2k-1}\tilde{\varepsilon}(t,r,\theta,q,\eta),
\end{gather}
the phase is $-\chi_0(\chi_1^2-\chi_2^{2k})$ which is the desired result. $\hfill{\blacksquare}$\medskip\\
If we use these local coordinates we obtain a simpler problem :
\begin{equation}
\int
e^{-\frac{i}{h}\chi_0(\chi_1^2-\chi_2^{2k})}A(\chi_0,\chi_1,\chi_2)
d\chi_0d\chi_1 d\chi_2,
\end{equation}
where the amplitude $A$ is obtained via pullback and integration,
i.e. :
\begin{equation}
A(\chi_0,\chi_1,\chi_2)=\int \chi^{*}(a |J\chi|)
d\chi_4...d\chi_{2n}.
\end{equation}
Now, we make several geometrical comments on this construction :
\begin{remark}\label{Jacobian}
\rm{The phase $-\chi_0(\chi_1^2-\chi_2^{2k})$ gives a correct
description of the singularities of the phase since the critical
set is :
\begin{equation}
\mathfrak{C}(-\chi_0(\chi_1^2-\chi_2^{2k}))=\{ (\chi_0,0,0)\}\cup %
\{ (0,\chi_1,\chi_1)\}.
\end{equation}
Moreover, this phase is non-degenerate w.r.t. $\chi_1$ when
$\chi_0\neq 0$. In order to express the main coefficient of the
trace formula, we notice that $\chi_2=0\Leftrightarrow x=x_0$ and
$\chi_1(t,0,\theta,q,\eta)=0 \Leftrightarrow \xi=0$, locally near
$(x_0,\xi_0)$. Also, because of the introduction of the polar
coordinates our amplitude satisfies :
\begin{equation*}
A(\chi_0,\chi_1,\chi_2)=\mathcal{O}(\chi_l ^{n-1}),\text{ }l=1,2.
\end{equation*}
Finally, each diffeomorphism used has Jacobian 1 in $z_0$ excepted
the correction w.r.t. $\alpha_2$, cf. Eq. (\ref{jacobian V_2K}),
which induces by pullback of the measure $r^{n-1}dr$ a
multiplication by :
\begin{equation}
|V_{2k}(\theta)|^{-\frac{n}{2k}}.
\end{equation}
These facts will be used in the next section to obtain an
invariant formulation.$\hfill{\square}$}
\end{remark}
\section{Proof of the main result.}
Let be $2k$ the even integer attached to our potential. We state 2
technical results which allow to compute the asymptotic expansion
of our oscillatory integrals. The first one proves the existence
of a total asymptotic expansion for the family of oscillatory
integrals attached to our normal forms. The second result, which
is the hardest part of this work, computes the main coefficients
of this expansion w.r.t. $h$ in the trace formula.
\begin{lemma} \label{technical result}
For $a\in C_0^\infty(\mathbb{R}\times [0,\infty]^2)$, the
oscillatory integrals :
\begin{equation}
J(\lambda)=\int\limits_{\mathbb{R}\times\mathbb{R}_{+}\times
\mathbb{R}_{+}} e^{-i\lambda \chi_0
(\chi_1^2-\chi_2^{2k})}A(\chi_0,\chi_1,\chi_2)d\chi_0 d\chi_1
d\chi_2,
\end{equation}
admit, as $\lambda\rightarrow \infty$, the asymptotic expansion :
\begin{equation} \label{structure DA}
J(\lambda) \sim \sum\limits_{j=0}^{\infty}
\lambda^{-\frac{l+1}{2k}} C_{j}(A)+\sum\limits_{j=0}^{\infty}
\lambda^{-j+1}\mathrm{log}(\lambda) D_{j}(A),
\end{equation}
where $C_{j}$ and $D_{j}$ are universal (computable)
distributions.
\end{lemma}
\begin{remark}\rm{The result stated obviously also holds for
integration on $\mathbb{R}\times\mathbb{R}^2$ if we split up the
domain of integration and use the symmetry w.r.t. $\chi_1$ and
$\chi_2$ of the phase. Also, note that terms with a logarithm of
the parameter only occur when $(j/2k)\in\mathbb{N}^*$.}
$\hfill{\square}$
\end{remark}
To compute the leading term of the trace formula we need a
particular and explicit result which explains the effect of the
dimension $n$ in our spectral problem. For $A\in
C_0^{\infty}(\mathbb{R}\times \mathbb{R}_{+}^2)$ we define :
\begin{equation}
I(\lambda)=\int\limits_{\mathbb{R}\times\mathbb{R}_{+}\times
\mathbb{R}_{+}} e^{-i\lambda t
(r^2-q^{2k})}A(t,r,q)r^{n-1}q^{n-1}dtdrdq.
\end{equation}
We note $t_{\pm}=\max (\pm,0)$ and $\hat{A}$ the partial Fourier
transform of $A$ w.r.t. $t$.
\begin{lemma} \label{Asympt. IO spectral problem}
If $n(k+1)/2k$ is not an integer then, as $\lambda\rightarrow
\infty$, we have :
\begin{equation} \label{IO E non-integer}
I(\lambda)= C_0(A)
\lambda^{-n\frac{k+1}{2k}}+\mathcal{O}(\log(\lambda)\lambda^{-\frac{n(k+1)+1}{2k}}).
\end{equation}
The distributional coefficients are given by :
\begin{gather*}
C_0(A)=C_{n,k}^{-} \int\limits_{\mathbb{R}}   t_{-}^{n\frac{k+1}{2k}-1}\hat{A}(t,0,0)dt,\text{ if $n$ is even,}\\
C_0(A)= \int\limits_{\mathbb{R}}  (C_{n,k}^{-}
t_{-}^{n\frac{k+1}{2k}-1}+C_{n,k}^{+}
t_{+}^{n\frac{k+1}{2k}-1})\hat{A}(t,0,0)dt, \text{ if $n$ is odd}.
\end{gather*}
But if $n(k+1)/2k\in\mathbb{N}^*$ and $n$ is odd, then :
\begin{equation}\label{IO E integer}
I(\lambda)=-C_0(A)
\lambda^{-n\frac{k+1}{2k}}\log(\lambda)+\mathcal{O}(\lambda^{-\frac{n(k+1)+1}{2k}}),
\end{equation}
with :
\begin{equation*}
C_0(A)=C_{n,k}\int\limits_{\mathbb{R}}
|t|^{n\frac{k+1}{2k}-1}\hat{A}(t,0,0)dt.
\end{equation*}
Finally, if  $n(k+1)/2k\in\mathbb{N}^*$ and $n$ is even the
asymptotic is given by Eq.(\ref{IO E non-integer}) with :
\begin{equation*}
C_0(A)= \int\limits_{\mathbb{R}}  (\tilde{C}_{n,k}^{-}
t_{-}^{n\frac{k+1}{2k}-1}+\tilde{C}_{n,k}^{+}
t_{+}^{n\frac{k+1}{2k}-1})\hat{A}(t,0,0)dt.
\end{equation*}
\end{lemma}
\begin{remark} \label{remark DA}\rm{We will not detail all the coefficients of the
asymptotic expansion because these are given by lengthy formulae
and it is, a priori, not possible to express invariantly their
contributions to the trace formula. Anyhow, these can be
explicitly computed with the procedure below. Also if
$\frac{n(k+1)+1}{2k} \notin \mathbb{N}$ the remainder of Eq.
(\ref{IO E non-integer}) can be optimized to
$\mathcal{O}(\lambda^{-\frac{n(k+1)+1}{2k}})$.} $\hfill{\square}$
\end{remark}
Before entering in the proof we would like to ad a comment
suggested by an interesting remark of D.Barlet. The Berstein-Sato
polynomial of our phase $\chi_0 (\chi_1^2 -\chi_2^{2k})$ can be
explicitly computed as :
\begin{equation}
\mathfrak{B}(z)=(z+1)^3 \prod\limits_{l=0, l\not= k+1}^{2k-2} (z +
\frac{l+k-1}{2k}).
\end{equation}
In particular, we could expect to obtain terms $\log(\lambda)^2$
in the expansion since $1/\mathfrak{B}$ has a triple pole in
$z=1$. But we will see that if the
amplitude is smooth there is no such contribution.\medskip\\
\textit{Proof of Lemma \ref{technical result}.} To attain our
objective we can restrict the proof of the lemma to an amplitude
$A(t,r,q)=a(t)B(r,q)$. This is justified by a standard density
argument in $C_{0}^{\infty}$ and the fact that the coefficients
obtained below are linear continuous functionals. By integration
w.r.t. $t$ we obtain :
\begin{equation*}
J(\lambda)=\int\limits_{\mathbb{R}_{+}^2} \hat{a}(\lambda
(r^2-q^{2k}))B(r,q)drdq.
\end{equation*}
This shows that the asymptotic is supported in the set $r=q^k$
since $\hat{a}$ decreases faster than any polynomial at
infinities. Although the new problem looks simple it is still too
complicated to obtain an explicit solution. First we split our
integral as :
\begin{equation}
J(\lambda)=J_{+}(\lambda)+J_{-}(\lambda),
\end{equation}
with :
\begin{gather*}
J_{+}(\lambda)=\int\limits_{0\leq q^k\leq r< \infty}
\hat{a}(\lambda
(r^2-q^{2k}))B(r,q)drdq,\\
J_{-}(\lambda)=\int\limits_{0\leq r\leq q^k <\infty}
\hat{a}(\lambda (r^2-q^{2k}))B(r,q)drdq.
\end{gather*}
Now we define the Melin transforms of $\hat{a}$ :
\begin{gather}
M_{+}(z)=\int\limits_{0}^{\infty} t^{z-1} \hat{a}(t)dt,\\
M_{-}(z)=\int\limits_{0}^{\infty} t^{z-1} \hat{a}(-t)dt.
\end{gather}
By Melin inversion formula we obtain that :
\begin{equation}
J_{+}(\lambda)=\frac{1}{2i\pi} \int\limits_{c+i\mathbb{R}}
M_{+}(z) \int\limits_{0\leq q^k\leq r< \infty} (\lambda
(r^2-q^{2k}))^{-z} B(r,q)drdq dz.
\end{equation}
Here $0<c<(2k)^{-1}$ so that the Melin inversion makes sense. In
order to desingularize the remaining part of the phase we
introduce the new coordinates $r=sq^k$, $s>0$ and $q>0$. We
accordingly obtain that :
\begin{equation}
J_{+}(\lambda)=\frac{1}{2i\pi} \int\limits_{c+i\mathbb{R}}
M_{+}(z) \lambda^{-z} \int\limits_{s=1}^{\infty}
\int\limits_{q=0}^{\infty} (s^2-1)^{-z} q^{-2kz} B(sq^k,q)dsq^kdq
dz,
\end{equation}
and :
\begin{equation}
J_{-}(\lambda)=\frac{1}{2i\pi} \int\limits_{c+i\mathbb{R}}
M_{-}(z) \lambda^{-z} \int\limits_{s=0}^{1}
\int\limits_{q=0}^{\infty} |s^2-1|^{-z} q^{-2kz} B(sq^k,q)dsq^kdq
dz.
\end{equation}
To control the remainders, and also to justify the changes of path
below, we recall the following classical result :
\begin{lemma}\label{Melin is S}
If $a\in \mathcal{S}(\mathbb{R})$ then we have :
\begin{equation}
\forall c>0 \text{ : } M_{\pm}(c+iy) \in\mathcal{S}(\mathbb{R}_y),
\end{equation}
a fortiori $M_{\pm}(c+iy)\in L^1(\mathbb{R},dy)$.
\end{lemma}
Our integrals can be expanded via Cauchy's residue method by
pushing of the complex path of integration to the right. The
associated distributional factors are analytic (cf. Lemma
\ref{distributions} below) and if we choose $d>c$ conveniently we
obtain :
\begin{equation*}
J_{+}(\lambda)= \sum\limits_{c<z_i<d} \mathrm{Res}
(z_i)+R(d,\lambda),
\end{equation*}
where the remainder satisfies :
\begin{equation}
|R(d,\lambda)|\leq C(B)\lambda^{-d} ||
M_{+}(d+iy)||_{L^1(\mathbb{R},dy)}=\mathcal{O}(\lambda^{-d}).
\end{equation}
Here, for each $d$ the constant $C(B)$ involves the $L^1$ norm of
a finite number of derivatives of $B$. This will indeed lead to an
asymptotic expansion.

The resulting asymptotics are related to poles of meromorphic
distributions $z\mapsto M_{\pm}(z) \lambda^{-z}(s^2-1)^{-z}
(q^{2k})^{-z}$ and it remains now to extend analytically these
distributions.
\begin{lemma} \label{distributions}
The family of distributions on $C_0^{\infty}(\mathbb{R}_{+}^2)$ :
$z\rightarrow |s^2-1|^{-z}q^{-2kz}$, initially defined in the
domain $\Re (z)<1/2k$, is meromorphic on $\mathbb{C}$ with poles
located at the rational points :
\begin{equation}z_{j,k}=\frac{j}{2k},\text{ } j\in\mathbb{N}^{*}.
\end{equation}
These poles are of order 2 if $z_{j,k}\in\mathbb{N}^*$ and of
order 1 otherwise.
\end{lemma}
\textit{Proof}. We first observe that for $s\geq 0$ we can write :
\begin{equation*}
|s^2-1|^{-z}=(s+1)^{-z} |s-1|^{-z}.
\end{equation*}
Since we will integrate on $s\geq0$ the first term defines clearly
an entire distribution. Now we remark that for $s\geq 1$ :
\begin{gather*}
\frac{\partial^{2k+1}}{\partial s \partial q^{2k}} (s-1)^{1-z}(q^{2k})^{1-z}=\mathfrak{b}_0(z)(s-1)^{-z}(q^{2k})^{-z},\\
\mathfrak{b}_0(z)=(1-z)\prod\limits_{j=1}^{2k} (j-2k).
\end{gather*}
Hence for $\Re(z)>0$ we can write :
\begin{gather} \notag
\int\limits_{s=1}^{\infty} \int\limits_{q=0}^{\infty}
(s^2-1)^{-z} q^{-2kz} f(s,q)dsdq \\
=\frac{-1}{\mathfrak{b}_0(z)} \int\limits_{s=1}^{\infty}
\int\limits_{q=0}^{\infty} (s-1)^{1-z}q^{2k(1-z)} \frac{\partial
^{2k+1}}{\partial s \partial q^{2k}} (1+s)^{-z}f(s,q)dsdq.
\end{gather}
Now, the r.h.s. is meromorphic in $\Re(z)<1+1/2k$ and we can
iterate to get the analytic continuation in $\Re(z)<l+1/2k$,
$\l\in\mathbb{N}$ arbitrary. The poles, and their order, can be
read off the rational functions of $z$ :
\begin{equation}
\mathfrak{R}_l(z)=\prod\limits_{m=0}^{l-1}
\frac{1}{\mathfrak{b}_0(z-m)}.
\end{equation}
Finally, a similar construction holds
if we integrate w.r.t. $s\in [0,1]$.$\hfill{\blacksquare}$\medskip\\
With the residue method and classical estimates for the remainder
(cf. Lemma \ref{Melin is S}) the last lemma proves the existence
of a total asymptotic expansion for the integrals $J(\lambda)$.
This ends the proof of Lemma \ref{technical result}.$\hfill{\blacksquare}$\medskip\\
\textit{Proof of Lemma \ref{Asympt. IO spectral problem}.} With
some slight modifications, we can apply the results above to our
initial problem. Taking Remark \ref{Jacobian} into account, we
have to compute asymptotics of oscillatory integrals with
amplitudes :
\begin{equation}
B(r,q)=b(r,q)r^{n-1}q^{n-1}.
\end{equation}
By substitution, we have to study the poles of :
\begin{gather}\label{analytic1}
\mathfrak{g}_{+}(z)=M_{+}(z)\int\limits_{s=1}^{\infty}
\int\limits_{q=0}^{\infty}
(s^2-1)^{-z}q^{-2kz} b(sq^k,q) s^{n-1} q^{n(k+1)-1} dsdq,\\
\mathfrak{g}_{-}(z)=M_{-}(z)\int\limits_{s=0}^{1}
\int\limits_{q=0}^{\infty} |s^2-1|^{-z}q^{-2kz} b(sq^k,q) s^{n-1}
q^{n(k+1)-1} dsdq. \label{analytic2}
\end{gather}
Here $b(sq^k,q)$ is in general no more of compact support but all
expressions and manipulations will be legal because terms
depending on $z$ will decrease faster and faster when we will
shift the path of integrations. To avoid unnecessary calculations
and discussions below, we remark that we can commute the
polynomial weights w.r.t. $q$ via the relation :
\begin{gather*}
\frac{\partial^{2k+1}}{\partial s
\partial q^{2k}}(s-1)^{1-z}q^{2k(1-z)}q^{n(k+1)-1}=
\mathfrak{b}(z) (s-1)^{-z}q^{-2kz} q^{n(k+1)-1},\\
\mathfrak{b}(z)=(z-1)\prod\limits_{j=1}^{2k} (j-2kz+n(k+1)-1).
\end{gather*}
More generally, with the differential operator:
\begin{equation*}
\mathrm{D}=\frac{\partial^{2k+1}}{\partial s
\partial q^{2k}},
\end{equation*}
after $l$-iterations one obtains :
\begin{gather*}
\mathrm{D}^l(s-1)^{l-z}q^{2k(l-z)}q^{n(k+1)-1}= \mathfrak{B}_l(z)
(s-1)^{-z}q^{-2kz} q^{n(k+1)-1},\\
\mathfrak{B}_l(z)=\prod\limits_{i=0}^{l-1}\mathfrak{b}(z-i),\text
{ } l\in\mathbb{N}^{*}.
\end{gather*}
A priori, this shows that there is poles of order 2 at positive
integers and simple poles at rational points :
\begin{equation}
z_{p,j,n,k}= p+\frac{j+n(k+1)-1}{2k}\notin \mathbb{N},\text{
}p\in\mathbb{N},\text{ }j\in [1,...,2k].
\end{equation}
For example, the analytic extension in the half-plane $\Re(z)<l$,
$l\in \mathbb{N}^{*}$, of the complex function appearing in the
first integral is explicitly given by :
\begin{equation} \label{full1}
\frac{(-1)^l}{\mathfrak{B}_l(z)}M_{+}(z)\lambda^{-z}\int\limits_{s=1}^{\infty}
\int\limits_{q=0}^{\infty} q^{2k(l-z)+n(k+1)-1} \mathrm{D}^l
(1+s)^{-z}b(sq^k,q)s^{n-1}dsdq.
\end{equation}
A similar relation for the term involving $M_{-}(z)$ is :
\begin{equation*}
\mathrm{D} |1-s|^{1-z} q^{n(k+1)-1-2k(1-z)}=-\mathfrak{b}(z)
|1-s|^{-z} q^{n(k+1)-1-2kz}.
\end{equation*}
This change of sign, due to the modulus, will be important below
because of some symmetries. After $l$-iterations we have :
\begin{equation}\label{full2}
\frac{1}{\mathfrak{B}_l(z)}M_{-}(z)\lambda^{-z}\int\limits_{s=0}^{1}
\int\limits_{q=0}^{\infty} q^{2k(l-z)+n(k+1)-1} \mathrm{D}^l
(1+s)^{-z}b(sq^k,q)s^{n-1}dsdq.
\end{equation}
\begin{remark}\label{roots b}\rm{The smallest double root of $\mathfrak{b}$,
and a fortiori of each $\mathfrak{B}_l$, is greater than
$z_{\rm{min}}=\frac{n(k+1)}{2k}$. In fact we will see below that
there is no poles, and a fortiori no contributions, before this
value. This insures that all integrals involved in the asymptotic
expansion are absolutely convergent.}$\hfill{\square}$
\end{remark}
A carefully examination of the integral w.r.t $q$ shows that all
coefficients are zero until we reach the pole :
\begin{equation}
z_{\mathrm{min}}=n\frac{k+1}{2k}=\frac{n}{2}+\frac{n}{2k}.
\end{equation}
This is justified by the fact that if $\alpha>0$, $\beta\in
\mathbb{N}^{*}$ and $\alpha+1>\beta$ we obtain :
\begin{equation*}
\int\limits_{0}^{\infty} \partial^\beta_x \left( x^\alpha f(x)
\right )dx =[\partial^{\beta-1}_x (x^\alpha
f(x))]_{x=0}^{\infty}=0,\text{ } \forall f\in
C_0^{\infty}(\mathbb{R}).
\end{equation*}
Since poles located at integers are of order 2 the attached
residuum are computed via the elementary formula :
\begin{equation}
\lim\limits_{z\rightarrow p} \frac{\partial}{\partial z}\left (
(z-p)^2 \mathfrak{g}_{\pm}(z)\lambda^{-z}\right ), \text{
}p\in\mathbb{N},
\end{equation}
where $\mathfrak{g}_{\pm}$ are defined by Eqs.
(\ref{analytic1},\ref{analytic2}). For $h$ holomorphic and
$\lambda>0$ we have :
\begin{equation}
\frac{\partial}{\partial z} (h(z)\lambda^{-z})=\frac{\partial
h}{\partial z}(z)\lambda^{-z}-\log(\lambda)\lambda^{-z}h(z),
\label{splitting}
\end{equation}
and we can apply this to $h(z)=(z-p)^2 \mathfrak{g}_{\pm}(z)$ near
$z=p$. Hence, a generic double pole located at $z=p$ leads to a
contribution :
\begin{equation}\label{coeff log^2}
(-1)^p \log(\lambda)\lambda^{-p}\lim\limits_{z\rightarrow p}
\frac{(z-p)^2}{\mathfrak{B}_p(z)}
M_{+}(z)\int\limits_{s=1}^{\infty} \int\limits_{q=0}^{\infty}
q^{n(k+1)-1}\mathrm{D}^p (s+1)^{-p} b(sq^k,q) s^{n-1}  dsdq.
\end{equation}
Since $M_{+}(p)$ is well defined and :
\begin{equation}
\lim\limits_{z\rightarrow
p}\frac{(z-p)^2}{\mathfrak{B}_p(z)}=c_p\in\mathbb{Q},
\end{equation}
our coefficients can be explicitly determined by the computation
of the integrals. Also, cf. Remark \ref{roots b}, the $c_p$ are
zero until $p\geq z_{\mathrm{min}}$. By integrations by parts and
up to a factorial number, the inner integral of Eq. (\ref{coeff
log^2}) can also be written :
\begin{gather*}
-\int\limits_{s=1}^{\infty} \frac{\partial^p}{\partial
s^p}(s+1)^{-p} s^{n-1}\left (\frac{\partial ^{2kp-1}}{\partial
q^{2kp-1}} b(sq^k,q)  q^{n(k+1)-1}  \right)_{|q=0} ds\\
=\left (\frac{\partial^{p-1}}{\partial s^{p-1}}(s+1)^{-p}
s^{n-1}\left (\frac{\partial ^{2kp-1}}{\partial q^{2kp-1}}
b(sq^k,q) q^{n(k+1)-1}  \right)_{|q=0}\right )_{|s=1}
\end{gather*}
But for all $g$ smooth bounded, with bounded derivatives, and $p$
large enough we have :
\begin{equation}
\int\limits_{0}^{\infty} \frac{\partial^p}{\partial s^p}\left(
(1+s)^{-p} s^{n-1} g(s)\right)ds=0.
\end{equation}
In particular, this is the case for all our coefficients since we
have :
\begin{equation*}
 p\geq\frac{n}{2}+\frac{n(k+1)}{2k}>n.
\end{equation*}
This trick shows that coefficients obtained by integration w.r.t.
$s$ on [0,1] and $[1,\infty]$ can be identified up to a sign and
we can save some computations.\medskip\\
\textbf{Computation of the leading term.}\\
Contrary to the case of non-degenerate critical points, the
evaluation of the leading term is somehow technical and some
computations will be left to the reader. According to the analysis
above, we distinguish out the case where the first non-zero
residue is attached to a simple or a
double pole.\medskip\\
\textbf{Case of $z_{\mathrm{min}}$ simple pole.}\\
Taking Remark \ref{roots b} into account the first coefficient is
given by :
\begin{gather*}
\lim\limits_{z\rightarrow z_{\mathrm{min}}}(-1)^l \lambda^{-z}
 \frac{(z-z_{\rm{min}})}{\mathfrak{B}_l(z)} M_{+}(z) \\
\int\limits_{s=1}^{\infty}
 \int\limits_{q=0}^{\infty} (s-1)^{l-z}q^{2k(l-z)+n(k+1)-1} \mathrm{D}^l (s+1)^{-z}b(sq^k,q)s^{n-1}
dsdq,
\end{gather*}
with $l\in \mathbb{N}^{*}$ such that $l>z_{\mathrm{min}}$ (any
such $l$ is acceptable). With this choice we can take the limit
under the integral to obtain :
\begin{gather}
c_{l,k,n} \lambda^{-n\frac{k+1}{2k}}
  \int\limits_{s=1}^{\infty}
 \int\limits_{q=0}^{\infty} (s-1)^{l-n\frac{k+1}{2k}}q^{2kl-1} \mathrm{D}^l (s+1)^{-n\frac{k+1}{2k}}b(sq^k,q)s^{n-1}
dsdq, \label{first coef}\\
c_{l,k,n}=(-1)^l M_{+}(n\frac{k+1}{2k})\lim\limits_{z\rightarrow
z_{\rm{min}}} \frac{(z-z_{\rm{min}})}{\mathfrak{B}_l(z)}.
\end{gather}
By integrations by parts, the integral w.r.t. $q$ of Eq.
(\ref{first coef}) is given by :
\begin{equation*}
\int\limits_{0}^{\infty} q^{2kl-1} \frac{\partial^{2kl}}{\partial
q^{2kl}} b(sq^k,q)dq=(2kl-1)! b(0,0).
\end{equation*}
Hence we obtain the asymptotic relation for the positive part of
our integral :
\begin{gather*}
I_{+}(\lambda)=c_{1}b(0,0)\lambda^{-n\frac{k+1}{2k}}
\int\limits_{1}^{\infty} (s-1)^{l-n\frac{k+1}{2k}}
\frac{\partial^l}{\partial s^l} ((1+s)^{-n\frac{k+1}{2k}}s^{n-1})
ds+R_1(\lambda),\\
c_{1}= (2kl-1)! c_{l,k,n}.
\end{gather*}
As concerns the convergence of the integral w.r.t. $s$, by
construction the singularity in $s=1$ is controlled.  For the
behavior at infinity we remark that the degree w.r.t. $s$ is :
\begin{equation}
l-n\frac{k+1}{2k}-n\frac{k+1}{2k}+n-1-l=-n\frac{k+1}{k} +n-1<-1,
\end{equation}
so that the integral is absolutely convergent. To compute
explicitly the values of our integrals we can choose $l=n$ since
$n-z_{\mathrm{min}}>0$. First, by induction on $n$ and for
$\frac{n}{2}<\Re(\alpha)<n+1$ we obtain that :
\begin{equation*}
E(n,\alpha)=\int\limits_{s=1}^{\infty} (s-1)^{n-\alpha}
\frac{\partial^n}{\partial s^n} ( (1+s)^{-\alpha} s^{n-1}) ds= 0,
\text{ if } n \text{ is even.}
\end{equation*}
Next, for $n$ odd, and always with $\frac{n}{2}<\Re(\alpha)<n+1$,
we have :
\begin{equation*} E(n,\alpha)=
 \prod\limits_{j=1}^{\frac{n-1}{2}}(-2j-1) 2^{\frac{n+1}{2}-2\alpha}\frac{\Gamma (n+1-\alpha) \Gamma(-n+2\alpha)}
{\Gamma(\frac{1-n}{2}+\alpha)}.
\end{equation*}
Hence for $\alpha=n(k+1)/2k$ and $n$ odd, we obtain :
\begin{gather}
\int\limits_{1}^{\infty} (s-1)^{n-n\frac{k+1}{2k}}
\frac{\partial^n}{\partial s^n} ((1+s)^{-n\frac{k+1}{2k}}s^{n-1})
ds=c_{n,k}\frac{\Gamma (1+n\frac{k-1}{2k})
\Gamma(\frac{n}{k})} {\Gamma(\frac{k+n}{2k})},\\
c_{n,k}=\prod\limits_{j=1}^{\frac{n-1}{2}}(-2j-1)
2^{\frac{n+1}{2}-n\frac{k+1}{k}}.
\end{gather}
For any integer $l>n(k+1)/2k$, similar computations show that :
\begin{gather*}
I_{-}(\lambda)=c_{2}b(0,0)\lambda^{-n\frac{k+1}{2k}}
\int\limits_{0}^{1} |s-1|^{l-n\frac{k+1}{2k}}
\frac{\partial^l}{\partial s^l} ((1+s)^{-n\frac{k+1}{2k}}s^{n-1})
ds+R_2(\lambda),\\
c_{2}= (2kl-1)! c_p M_{-}(\frac{n(k+1)}{2k}).
\end{gather*}
Once more the choice of $l=n$ is admissible and we define :
\begin{equation*}
a_{n,k}=\int\limits_{0}^{1} |s-1|^{n-n\frac{k+1}{2k}}
\frac{\partial^n}{\partial s^n} ((1+s)^{-n\frac{k+1}{2k}}s^{n-1})
ds.
\end{equation*}
since these integrals seem to have no formulation by mean of
elementary functions, unless by mean of hypergeometric functions.
These numbers $a_{n,k}$ are finite and non-zero in general
position.

From the analysis of the poles above we know that the remainders
$R_1$ and $R_2$ are of order
$\mathcal{O}(\lambda^{-\frac{n(k+1)+1}{2k}})$ if the next pole is
not an integer, respectively $\mathcal{O}(\log(\lambda)
\lambda^{-\frac{n(k+1)+1}{2k}})$ if this is an integer (cf. Remark
\ref{remark DA}). By summation we obtain the leading term of the
asymptotic expansion with a precise remainder.
\medskip\\
\textbf{Case of $z_{\mathrm{min}}$ pole of order 2.}\\
Starting from Eq. (\ref{splitting}) we see that the associated
coefficients are given by :
\begin{equation}
-\log(\lambda) \lambda^{-z_{\mathrm{min}}}
\left((z-z_{\mathrm{min}})^2 \mathfrak{g}_{\pm} (z) M_{\pm}(z)
\right) _{|z=z_{\mathrm{min}}}.
\end{equation}
But a great part of this limit was precisely computed above. Since
$z_{\mathrm{min}}$ is by assumption an integer we will obtain some
particular values. By induction on $p>1$, and assuming recursively
that $n<2p$, we have :
\begin{equation}
\int\limits_{1}^{\infty} \partial_s^{p} ((s+1)^{-p}
s^{n-1})ds=-\frac{1}{2^p}\prod\limits_{j=0}^{p-1} (n-2j).
\end{equation}
If $n$ is odd this coefficient is not zero and we get the result.
But if $n$ is even we obtain that the associated contribution
vanishes and there is no logarithm in the leading term. To obtain
the top order coefficient we must compute the coefficient obtained
by derivation of our meromorphic distributions.

Starting from Eqs. (\ref{full1},\ref{full2}), by Leibnitz rule, we
have 3 possibilities : derivation of the rational function, of the
Melin transform or of the analytic integrals in $(q,s)$. Since the
integral w.r.t. $s$ vanishes in $z=z_{\mathrm{min}}$, the 2 first
terms do not contribute. Similarly, the derivative of $q^{-2kz}$
can be discarded. Hence, the only contribution comes from
derivation of the distribution w.r.t. $s$ and we have to use the
modified constants :
\begin{gather}
\tilde{a}^{+}_{p,k}=-\int\limits_{1}^{\infty} \log(s^2-1)
\partial_s^{p} ((s+1)^{-p}s^{n-1}) ds,\\
\tilde{a}^{-}_{p,k}=-\int\limits_{0}^{1} \log(s^2+1)
\partial_s^{p} ((s+1)^{-p}s^{n-1}) ds,
\end{gather}
respectively for $I_{+}(\lambda)$ and $I_{-}(\lambda)$. By similar
considerations as above, these integrals are absolutely convergent for any $p\geq z_{\mathrm{min}}$.\medskip\\
\textbf{The other coefficients.}\\
To obtain a complete overview of the asymptotic expansion we show
also how to compute the coefficients attached to logarithmic
distributions. For $p\in\mathbb{N}^{*}$, $p\geq z_{\mathrm{min}}$,
by Leibnitz rule, these derivatives are equal to :
\begin{equation*}
\partial_z M_{\pm}(p) \lim\limits_{z\rightarrow p} (z-p)^2\mathfrak{g}_{\pm}(z)+
M_{\pm}(p) \lim\limits_{z\rightarrow p}
(\partial_z(z-p)^2\mathfrak{g}_{\pm}(z)).
\end{equation*}
Hence the asymptotic involves also the distributional terms :
\begin{equation*}
\int\limits_{0}^{\infty} \log(t) t^{p-1} \hat{a}(\pm t)dt, \text{
}p\in\mathbb{N}^{*},
\end{equation*}
and, by derivation of the meromorphic distributions, terms :
\begin{gather*}
\int\limits_{s=1}^{\infty} \int\limits_{q=0}^{\infty} \log (q)
q^{\alpha} D^{p}(1+s)^{-p}b(sq^k,q) s^{n-1} dsdq,\\
\int\limits_{s=1}^{\infty} \int\limits_{q=0}^{\infty} \log(s^2-1)
q^{\alpha} D^{p}(1+s)^{-p}b(sq^k,q) s^{n-1} dsdq.
\end{gather*}
Where the parameter $\alpha$ runs in the sequence of positive
rational numbers $l/2k$. There is also similar terms with
integration w.r.t. $s\in [0,1]$. Also, from the analysis above, we
know that the logarithmic coefficients only occur when
$p=l(k+1)/2k$ are integers.\medskip\\
\textbf{Invariant formulation of the main coefficients.}\\
To complete the proof it remains to express our coefficients in
terms of the natural data of the problem, i.e. $\varphi$ and $V$.
We apply Lemma \ref{technical result} to our amplitude
$A(\chi_0,\chi_1,\chi_2)$ to prove the existence of the total
asymptotic expansion. As concerns the leading term, in any of the
cases at hand we have to express the amplitude located at the
origin. By standard manipulations, already used in
\cite{Cam1,Cam2}, we can inverse our diffeomorphism via an
oscillatory representation of the delta-Dirac distribution by mean
of a Schwartz kernel :
\begin{equation*}
K(\delta_{\{\chi_1,\chi_2\}})=\frac{1}{(2\pi)^2}
\int\limits_{\mathbb{R}^2} e^{-i \chi_1 z_1} e^{-i \chi_2 z_2}
dz_1 dz_2.
\end{equation*}
After integration w.r.t. $d\theta$, we accordingly obtain that :
\begin{equation}\label{ampli}
b(t,0,0)=\mathrm{S}(\mathbb{S}^{n-1}) a(t,z_0)
\int\limits_{\mathbb{S}^{n-1}} |V_{2k}(\eta)|^{-\frac{n}{2k}}
d\eta.
\end{equation}
We recall that these integrals over the spheres are simply given
by the Jacobian of our coordinates on the blow-up of the critical
point (cf. Remark \ref{Jacobian}). The principal symbol of
$\Theta(P_h)$ is $\Theta(p)$ with $\Theta(p(z_0))=1$. Hence, at
the critical point, the term homogeneous of degree 0 w.r.t. $h$ of
the amplitude of our FIO is given by $a(t,z_0)=\hat{\varphi}(t)$
(cf. section 3). Substituting Eq. (\ref{ampli}) in all integral
formulas for the leading terms of the asymptotic expansion the
Fourier inversion formula yields :
\begin{gather*}
M_{+}(\hat{a}(t,z_0))(n\frac{k+1}{2k})=\int\limits_{0}^{\infty} \varphi(t) t^{n\frac{k+1}{2k}-1}dt,\\
M_{-}(\hat{a}(t,z_0))(n\frac{k+1}{2k})=\int\limits_{0}^{\infty}
\varphi(-t) t^{n\frac{k+1}{2k}-1}dt.
\end{gather*}
Setting $\lambda=h^{-1}$, so that $\log(h)=-\log(\lambda)$,
dividing by $(2\pi h)^n$ we obtain, via Lemma \ref{Asympt. IO
spectral problem}, the
results stated in Theorem \ref{Main}. $\hfill{\blacksquare}$\medskip\\
\textbf{Extension.}\medskip\\
We show here shortly how to extend the result of Theorem
\ref{Main} to the case of an $h$-admissible operator $P_h$ of
symbol $p_h\sim \sum h^j p_j$ whose principal symbol is
$p_0=\xi^2+V(x)$ with a non-vanishing subprincipal symbol $p_1$.
In this case the Fourier integral operator approximating the
propagator has the amplitude :
\begin{equation*}
\tilde{a}(t,z)= a(t,z)\exp (i \int\limits_{0}^{t}
p_1(\Phi_s(z))ds),
\end{equation*}
see Duistermaat \cite{DUI1} concerning the solution of the first
transport equation. Since $z_0$ is an equilibrium we have simply
$p_1(\Phi_s(z_0))=p_1(z_0)$ and hence :
\begin{equation}
\tilde{a}(t,z_0)= \hat{\varphi}(t) e^{it p_1(z_0)}.
\end{equation}
If the subprincipal symbol vanishes at the critical point, which
is the case in a lot of practical situations, the top order
coefficient of the trace formula remains the same. Finally, when
$p_1(z_0)\neq 0$ by Fourier inversion formula we replace
$\varphi(t)$ by $\varphi(t+p_1(z_0))$ in all integral formulae of
Theorem \ref{Main}.

\end{document}